%% file: main.tex
\documentclass[11pt]{amsart}
\usepackage{amssymb}
\usepackage[colorlinks=true,linkcolor=blue,citecolor=blue,urlcolor=blue]{hyperref}
\newtheorem{theorem}{Theorem}[section]
\newtheorem{lemma}[theorem]{Lemma}
\newtheorem{proposition}[theorem]{Proposition}
\newtheorem{corollary}[theorem]{Corollary}
\newtheorem*{liuconj}{Conjecture}
\theoremstyle{remark}
\newtheorem{remark}[theorem]{Remark}

\newcommand{\F}{\mathbb{F}}
\newcommand{\Tr}{\operatorname{Tr}}
\newcommand{\PG}{\mathrm{PG}}

\input{tables}

\begin{document}

\title[{Two conjectures of Liu}]{Two conjectures of Liu on self-orthogonal MDS codes}
\author{Vadim Genkin}
\address{Independent researcher, Tallinn, Estonia}
\email{vadim.genkin@gmail.com}
\date{26 November 2026}
\subjclass[2020]{94B05; 51E21, 05B25}
\keywords{MDS code, self-orthogonal code, hyperoval, translation hyperoval, Reed--Solomon code, arc, normal rational curve}

\begin{abstract}
Shu Liu recently studied for which parameters $q$-ary self-orthogonal MDS codes exist, and made two conjectures. Conjecture~1 states that there is no $[q+2,3]$ self-orthogonal MDS code for even $q\ge4$. Assuming the Main MDS Conjecture, this is the only case left open in his classification for even $q$, and he proved it for $q=4$. We show that Conjecture~1 also holds for $q=8$ and fails for every even $q\ge16$. For $q=2^m$ with $m=5$ or $m\ge7$ the counterexamples form an explicit family: if $2\le i\le m-2$, $\gcd(i,m)=1$ and $c\in\F_q$ has absolute trace~$1$, then the matrix with columns $(c+t^{-1}+t^{-2^i})\,(1,t,t^{2^i})^{T}$ for $t\in\F_q$ (where $0^{-1}=0$), $(0,1,0)^{T}$ and $(0,0,1)^{T}$ generates a self-orthogonal $[q+2,3,q]$ MDS code. For $q=16$ and $q=64$ we give explicit codes. Conjecture~2 states that there is no $[q+1,(q-1)/2]$ self-orthogonal MDS code for odd $q$; Liu proved it for $q=3$ and $q=5$. We prove it for every prime $q$ and for $q=9$. By theorems of Ball and of Segre, every such MDS code is then, up to the order of the coordinates, a doubly extended generalized Reed--Solomon code, and we show that no such Reed--Solomon code is self-orthogonal, for any odd $q$. So for the remaining odd $q$, Conjecture~2 would follow from the conjecture that all arcs of size $q+1$ in these dimensions are normal rational curves.
\end{abstract}

\maketitle

\section{Introduction}

A linear $[n,k,d]$ code $C\subseteq\F_q^n$ is \emph{MDS} if $d=n-k+1$, and (Euclidean) \emph{self-orthogonal} if $C\subseteq C^\perp$, where $C^\perp=\{v\in\F_q^n:\sum_ju_jv_j=0\text{ for all }u\in C\}$. If $G$ is a generator matrix of $C$, then $C$ is self-orthogonal if and only if $GG^{T}=0$.

Shu Liu~\cite{Liu} studies for which $q$, $n$ and $k$ a $q$-ary $[n,k]$ self-orthogonal MDS code exists, and makes two conjectures. For even $q$, assuming the Main MDS Conjecture, he solves this problem except for the case $[n,k]=[q+2,3]$ \cite[\S1.2 and Main Theorem~1]{Liu}, proves that no $4$-ary $[6,3]$ self-orthogonal MDS code exists \cite[Proposition~4.7]{Liu}, and conjectures the following.

\begin{liuconj}[{\cite[Conjecture~1]{Liu}}]
For an even $q\geqslant4$, there do not exist $q$-ary $[q+2,3]$ self-orthogonal MDS codes.
\end{liuconj}

For even $q$ we write $q=2^m$ with $m\ge1$, and $\Tr\colon\F_q\to\F_2$ is the absolute trace. A $[q+2,3]$ MDS code has minimum distance $q$, and no $3$-dimensional MDS code over $\F_q$ is longer (Lemma~\ref{lem:mds}).

\begin{theorem}\label{thm:main}
Let $q\ge4$ be even. A $q$-ary $[q+2,3]$ self-orthogonal MDS code exists if and only if $q\ge16$. Hence Liu's Conjecture~1 is true for $q=4$ \textup{(}as proved in \cite[Proposition~4.7]{Liu}\textup{)} and $q=8$, and false for every even $q\ge16$.
\end{theorem}

Together with \cite{Liu}, this completes the classification for even $q$, under the Main MDS Conjecture as assumed in \cite{Liu}. Independently, and slightly earlier, Moreno-Mej\'{\i}a and Jim\'enez-Hern\'andez made public explicit self-orthogonal $[q+2,3]$ MDS codes for $q=16$, $32$, $64$ and $128$ \cite{MJ}.

The main part of the proof is an explicit family.

\begin{theorem}\label{thm:family}
Let $q=2^m$, let $i$ be an integer with $2\le i\le m-2$ and $\gcd(i,m)=1$, put $E=2^i$, and choose $c\in\F_q$ with $\Tr(c)=1$. For $t\in\F_q$ let
\[
\lambda(t)=c+t^{-1}+t^{-E},\qquad\text{with the convention }0^{-1}=0\text{ (so }\lambda(0)=c\text{)}.
\]
Let $G$ be the $3\times(q+2)$ matrix whose columns are $\lambda(t)\,(1,t,t^{E})^{T}$ for $t\in\F_q$, followed by $(0,1,0)^{T}$ and $(0,0,1)^{T}$. Then the row space of $G$ is a self-orthogonal $[q+2,3,q]$ MDS code over $\F_q$. An integer $i$ as above exists if and only if $m=5$ or $m\ge7$.
\end{theorem}

The columns of $G$ represent the points of the translation hyperoval $D(2^i)$. The columns of a self-orthogonal $[q+2,3]$ MDS code never represent a regular hyperoval (Corollary~\ref{cor:regular}). For $m\in\{4,6\}$ every hyperoval $D(2^i)$ with $\gcd(i,m)=1$ is regular (Remark~\ref{rem:boundary}), and Theorem~\ref{thm:family} does not apply. For these two fields we give explicit matrices (Section~\ref{sec:1664}); their columns represent the Lunelli--Sce hyperoval for $q=16$ and a Subiaco hyperoval for $q=64$. For $q\in\{4,8\}$ every hyperoval is regular, and no code exists (Section~\ref{sec:small}).

Up to the order of the coordinates, the $[q+2,3]$ MDS codes whose columns represent a regular hyperoval are the generalized triply extended Reed--Solomon codes, obtained from a $3$-dimensional doubly extended generalized Reed--Solomon code by appending a nonzero multiple of the column $(0,1,0)^{T}$, the nucleus of the conic. By Corollary~\ref{cor:regular} none of them is self-orthogonal, so the codes of Theorem~\ref{thm:main} are not of this type.

For odd $q$, Liu settles the existence of $q$-ary $[n,k]$ self-orthogonal MDS codes with $k\le n/2$ and $n=q+1$, $n\mid q$ or $n\mid(q-1)$, except in two cases \cite[Main Theorem~2 and Remark~1.1]{Liu}. For the first of them he conjectures the following. He proves it for $q=3$ and $q=5$ \cite[Proposition~4.9]{Liu}, and, assuming the Main MDS Conjecture, for codes contained in a self-dual MDS code \cite[Proposition~4.10]{Liu}.

\begin{liuconj}[{\cite[Conjecture~2]{Liu}}]
For odd $q$, there do not exist any $[n,k]=[q+1,(q-1)/2]$ self-orthogonal MDS codes.
\end{liuconj}

\begin{theorem}\label{thm:odd}
Let $q$ be odd. No doubly extended generalized Reed--Solomon code of length $q+1$ and dimension $(q-1)/2$ is self-orthogonal. Consequently, Liu's Conjecture~2 holds for $q$ if $q$ is a prime or $q=9$, and, for $q\ge5$, whenever every arc of size $q+1$ in $\PG((q-3)/2,q)$ is a normal rational curve.
\end{theorem}

The proof is in Section~\ref{sec:odd}. A short computation reduces the first statement to the fact that a quadratic polynomial over $\F_q$ cannot take only nonzero square values. For prime $q$ and for $q=9$, theorems of Ball~\cite{Ba} and Segre~\cite{Se55} show that, up to the order of the coordinates, every $[q+1,(q-1)/2]$ MDS code is a doubly extended generalized Reed--Solomon code. Ball and Lavrauw conjecture that every arc of size $q+1$ in $\PG(k-1,q)$ with $6\le k\le q-5$ is a normal rational curve \cite[Conjecture~22]{BL}. For $q\ge13$ the dimension $k=(q-1)/2$ is in this range, so their conjecture would give Conjecture~2 for every odd $q$, much as Liu's classification for even $q$ relies on the Main MDS Conjecture. We do not treat the second case left open in \cite[Main Theorem~2]{Liu}, namely $[n,k]=[q-1,(q-1)/2]$ with $q\equiv1\pmod4$.

\section{Preliminaries}\label{sec:prelim}

In Sections~\ref{sec:prelim}--\ref{sec:small}, $\F=\F_q$ with $q=2^m$. A point of the projective plane $\PG(2,q)$ is a $1$-dimensional subspace $\langle v\rangle$ of $\F^3$; its \emph{normalized} representative is the one whose first nonzero coordinate is $1$. An \emph{arc} is a set of points no three of which are collinear, and a \emph{hyperoval} is an arc of $q+2$ points \cite{Hir}.

\begin{lemma}\label{lem:mds}
Let $G=[g_1|\cdots|g_n]$ be a $3\times n$ matrix over $\F$ with $n\ge3$. The row space of $G$ is an $[n,3,n-2]$ MDS code if and only if every three columns of $G$ are linearly independent. In that case the points $\langle g_1\rangle,\dots,\langle g_n\rangle$ are distinct and form an arc, so $n\le q+2$.
\end{lemma}

\begin{proof}
If every three columns are independent, $G$ has rank $3$. For $x\ne0$, the zero coordinates of the codeword $xG$ are the indices $j$ with $g_j$ in the $2$-dimensional subspace $x^\perp=\{y:\sum_ax_ay_a=0\}$, which contains at most two columns; so $d\ge n-2$, and $d=n-2$ by the Singleton bound. Conversely, if three columns lie in a subspace $x^\perp$ with $x\ne0$, then $xG$ is a nonzero codeword with at least three zero coordinates. No two columns are proportional, and no three represent collinear points. Fixing one point $P$ of the arc, each of the other $n-1$ points lies on one of the $q+1$ lines through $P$, and each line contains at most one of them; hence $n\le q+2$.
\end{proof}

A \emph{quadratic form} is $Q(x_1,x_2,x_3)=\sum_{1\le a\le b\le3}Q_{ab}x_ax_b$ with $Q_{ab}\in\F$; it satisfies $Q(\mu v)=\mu^2Q(v)$.

\begin{lemma}[Gram identity]\label{lem:gram}
For every $3\times n$ matrix $G=[g_1|\cdots|g_n]$ and every quadratic form $Q=\sum_{a\le b}Q_{ab}x_ax_b$, $\sum_{j}Q(g_j)=\sum_{a\le b}Q_{ab}(GG^{T})_{ab}$. In particular, $GG^{T}=0$ implies $\sum_jQ(g_j)=0$.
\end{lemma}

\begin{proof}
$(GG^{T})_{ab}=\sum_jg_j[a]\,g_j[b]$; multiply by $Q_{ab}$ and sum over $a\le b$.
\end{proof}

\begin{lemma}[scalings]\label{lem:search}
Let $v_1,\dots,v_n$ represent the points of a hyperoval, $n=q+2$, and let $\lambda\in(\F^*)^n$. Then $G_\lambda=[\lambda_1v_1|\cdots|\lambda_nv_n]$ generates an MDS code, and $G_\lambda G_\lambda^{T}=\sum_j\lambda_j^2v_jv_j^{T}$. So $G_\lambda$ generates a self-orthogonal MDS code if and only if $w=(\lambda_j^2)_j$ satisfies the six linear equations $\sum_jw_jv_j[a]v_j[b]=0$, $1\le a\le b\le3$. Since squaring is a bijection of $\F$, every solution $w$ with all $w_j\neq0$ gives such a code, with $\lambda_j=w_j^{q/2}$.
\end{lemma}

\begin{proof}
Every $3\times3$ minor of $G_\lambda$ is a product of three $\lambda_j$ and a nonzero determinant; apply Lemma~\ref{lem:mds}. The formula for $G_\lambda G_\lambda^{T}$ is immediate.
\end{proof}

\begin{lemma}[nucleus obstruction]\label{lem:nucleus}
Let $v_1,\dots,v_n\in\F^3$ and let $Q$ be a quadratic form with $Q(v_j)=0$ for all $j\ne N$ and $Q(v_N)\ne0$. If $G=[\lambda_1v_1|\cdots|\lambda_nv_n]$ satisfies $GG^{T}=0$, then $\lambda_N=0$.
\end{lemma}

\begin{proof}
By Lemma~\ref{lem:gram}, $0=\sum_jQ(\lambda_jv_j)=\sum_j\lambda_j^2Q(v_j)=\lambda_N^2Q(v_N)$.
\end{proof}

A hyperoval is \emph{regular} if it is projectively equivalent to $H_{\rm reg}=\{(1,t,t^2):t\in\F\}\cup\{(0,1,0),(0,0,1)\}$, a conic together with its nucleus $(0,1,0)$ \cite{Hir}. The form $Q_0=x_2^2+x_1x_3$ vanishes at every point of $H_{\rm reg}$ except the nucleus, where it equals $1$. For $A\in GL(3,q)$, $Q_0\circ A^{-1}$ is a quadratic form with $(Q_0\circ A^{-1})(Av)=Q_0(v)$, so it vanishes at every point of $AH_{\rm reg}$ except $\langle A(0,1,0)^{T}\rangle$. With Lemmas~\ref{lem:mds} and~\ref{lem:nucleus} this gives:

\begin{corollary}\label{cor:regular}
The columns of a generator matrix of a self-orthogonal $[q+2,3]$ MDS code never represent the points of a regular hyperoval.
\end{corollary}

\section{The explicit family}\label{sec:family}

For $i\ge1$ let $D(2^i)=\{(1,t,t^{2^i}):t\in\F\}\cup\{(0,1,0),(0,0,1)\}$.

\begin{lemma}[Segre~\cite{Se57}]\label{lem:translation}
If $\gcd(i,m)=1$, then any three vectors of $D(2^i)$ are linearly independent, so $D(2^i)$ represents a hyperoval.
\end{lemma}

\begin{proof}
Put $E=2^i$. Since $\gcd(2^i-1,2^m-1)=2^{\gcd(i,m)}-1=1$, the map $x\mapsto x^{E-1}$ is injective on $\F^*$. For distinct $a,b,a'\in\F$, put $u=a+b$ and $v=a+a'$. By additivity of $t\mapsto t^E$, the determinant of $(1,a,a^E),(1,b,b^E),(1,a',a'^E)$ is $uv^E+u^Ev=uv(u^{E-1}+v^{E-1})$, which is nonzero since $u,v\ne0$, $u\ne v$ and $x\mapsto x^{E-1}$ is injective. The determinants of $(1,a,a^E),(1,b,b^E),(0,1,0)$, of $(1,a,a^E),(1,b,b^E),(0,0,1)$ and of $(1,a,a^E),(0,1,0),(0,0,1)$ are $(a+b)^E$, $a+b$ and $1$.
\end{proof}

\begin{lemma}[power sums]\label{lem:powersum}
For $r\in\mathbb Z$ let $S(r)=\sum_{t\in\F^*}t^r$. Then $S(r)=1$ if $(q-1)\mid r$, and $S(r)=0$ otherwise.
\end{lemma}

\begin{proof}
If $(q-1)\mid r$, all $q-1$ terms are $1$ and $q-1$ is odd. Otherwise, for a generator $\gamma$ of $\F^*$ we have $\gamma^r\ne1$ and $\gamma^rS(r)=S(r)$.
\end{proof}

\begin{proof}[Proof of Theorem~\ref{thm:family}]
\emph{$c$ exists.} $\Tr$ is a nonzero polynomial of degree $2^{m-1}<q$ with values in $\F_2$, so $\Tr(c)=1$ for some $c$.

\emph{$\lambda$ never vanishes.} For $t\ne0$ put $s=t^{-1}$; then $t^{-E}=s^{2^i}$ and $\Tr(s^{2^i})=\Tr(s)$, so $\Tr(\lambda(t))=\Tr(c)+\Tr(s)+\Tr(s^{2^i})=1$. Also $\Tr(\lambda(0))=\Tr(c)=1$. Hence $\lambda(t)\ne0$ for all $t$.

\emph{MDS.} By Lemma~\ref{lem:translation} and Lemma~\ref{lem:search}, the row space of $G$ is a $[q+2,3,q]$ MDS code.

\emph{$GG^{T}=0$.} For $e\ge0$ let $T(e)=\sum_{t\in\F}\lambda(t)^2t^e$, with $0^0=1$. Then
\[
GG^{T}=\begin{pmatrix}T(0)&T(1)&T(E)\\T(1)&T(2)+1&T(E+1)\\T(E)&T(E+1)&T(2E)+1\end{pmatrix}.
\]
Since $\lambda(t)^2=c^2+t^{-2}+t^{-2E}$ for $t\ne0$ and $\lambda(0)^2=c^2$, we get
\[
T(e)=c^2[e=0]+c^2S(e)+S(e-2)+S(e-2E),
\]
where $[P]$ is $1$ if $P$ holds and $0$ otherwise. As $2\le i\le m-2$, we have $4\le E\le q/4$, so for $e\in\{0,1,2,E,E+1,2E\}$ the arguments $e$, $e-2$, $e-2E$ lie in $[-q/2,q/2]$, where the only multiple of $q-1$ is $0$ (as $q-1>q/2$). By Lemma~\ref{lem:powersum},
\[
T(e)=c^2[e=0]+c^2[e=0]+[e=2]+[e=2E]=[e=2]+[e=2E].
\]
The numbers $0,1,2,E,E+1,2E$ are distinct, so $T(2)=T(2E)=1$ and the other four values are $0$; hence $GG^{T}=0$. (This step does not use $\gcd(i,m)=1$, which is needed only for the MDS property.)

\emph{Admissible $i$.} For odd $m\ge5$ take $i=2$. For $m\equiv0\pmod4$, $m\ge8$, take $i=m/2-1$, which is odd; a common divisor of $i$ and $m$ divides $m-2i=2$ and is odd, so it is $1$. For $m\equiv2\pmod4$, $m\ge10$, take $i=m/2-2$, which is odd; a common divisor divides $m-2i=4$ and is odd, so it is $1$. In each case $2\le i\le m-2$. Conversely, for $m\le3$ the range $2\le i\le m-2$ is empty, and for $m=4$ and $m=6$ every $i$ in the range has $\gcd(i,m)>1$.
\end{proof}

\begin{remark}\label{rem:boundary}
For $i=1$ and $i=m-1$ the set $D(2^i)$ is regular: $D(2)=H_{\rm reg}$, and $D(2^{m-1})$ is $H_{\rm reg}$ with the last two coordinates swapped (substitute $t=s^2$). For $m\in\{4,6\}$, every $i$ with $1\le i\le m-1$ and $\gcd(i,m)=1$ is $1$ or $m-1$, so every hyperoval $D(2^i)$ with $\gcd(i,m)=1$ is regular and gives no code.
\end{remark}

\section{The fields with 16 and 64 elements}\label{sec:1664}

We write $\F_{2^m}=\F_2[x]/(p)$ and encode $\sum_kb_kx^k$ ($b_k\in\{0,1\}$) as the integer $\sum_kb_k2^k$; we use $p=x^4+x+1$ for $q=16$ and $p=x^6+x+1$ for $q=64$, and $\eta$ denotes the class of $x$.

\begin{proposition}\label{prop:1664}
The following matrices $G_{16}$ over $\F_{16}$ and $G_{64}$ over $\F_{64}$ generate self-orthogonal $[18,3,16]$ and $[66,3,64]$ MDS codes:
\[
\setlength{\arraycolsep}{2.5pt}
G_{16}=\MatrixSixteen,
\]
and $G_{64}=[A\,|\,B\,|\,C]$ with
{\footnotesize
\[
\setlength{\arraycolsep}{1.6pt}
A=\MatrixSixtyFourA,
\]
\[
\setlength{\arraycolsep}{1.6pt}
B=\MatrixSixtyFourB,
\]
\[
\setlength{\arraycolsep}{1.6pt}
C=\MatrixSixtyFourC.
\]}
\end{proposition}

\begin{proof}
By Lemma~\ref{lem:mds}, this amounts to $GG^{T}=0$ and the nonvanishing of the $\binom{18}{3}=816$ and $\binom{66}{3}=45760$ column-minors, which we verified by computer (Section~\ref{sec:comp}).
\end{proof}

The normalized columns of $G_{16}$ are the points $(1,t,f(t))$, $t\in\F_{16}$, in the order of the integer encoding, followed by $(0,1,0)$ and $(0,0,1)$, where $f(t)=t^{12}+t^{10}+\eta^{11}t^8+t^6+\eta^2t^4+\eta^9t^2$ is an o-polynomial of the Lunelli--Sce hyperoval \cite{LS}. Likewise, the normalized columns of $G_{64}$ are the points $(1,t,f(t))$, $t\in\F_{64}$, in the order of the integer encoding, followed by $(0,1,0)$ and $(0,0,1)$, where
\[
f(t)=\frac{d^2(t^4+t)+d^2(1+d+d^2)(t^3+t^2)}{t^4+d^2t^2+1}+t^{1/2},\qquad d=\eta,
\]
is an o-polynomial of a Subiaco hyperoval \cite{CPPR}. Here $d\notin\F_4$ and $\Tr(1/d)=1$, as required in \cite{CPPR}, and the denominator $(t^2+\eta t+1)^2$ has no zero in $\F_{64}$. These identifications explain where the matrices come from, but the proof does not use them. We found both matrices with Lemma~\ref{lem:search}: we computed the solution space of the six equations for these hyperovals and drew random solutions until one had no zero coordinate.

\section{The fields with at most 8 elements}\label{sec:small}

For $q=2$ ($m=1$) there is no code, since a self-orthogonal code has $k\le n/2$. The \emph{standard frame} consists of the points $(1,0,0)$, $(0,1,0)$, $(0,0,1)$ and $(1,1,1)$.

\begin{proposition}\label{prop:nonexist}
Suppose that every hyperoval $H$ of $\PG(2,q)$ that contains the standard frame admits a quadratic form $Q_H$ and a point $N_H\in H$ such that $Q_H$ vanishes at every point of $H$ other than $N_H$, but not at $N_H$. Then no $q$-ary $[q+2,3]$ self-orthogonal MDS code exists.
\end{proposition}

\begin{proof}
Let $G=[g_1|\cdots|g_n]$ generate such a code. By Lemma~\ref{lem:mds}, $g_4=\alpha_1g_1+\alpha_2g_2+\alpha_3g_3$ with $\alpha_1,\alpha_2,\alpha_3\ne0$. Let $A$ be the inverse of $[\alpha_1g_1|\alpha_2g_2|\alpha_3g_3]$ and $G'=AG$. Then $G'$ has the same row space as $G$, $G'G'^{T}=A(GG^{T})A^{T}=0$, and the first four columns of $G'$ represent the standard frame. By Lemma~\ref{lem:mds} the columns of $G'$ represent a hyperoval $H$ through the frame. Writing them as $\lambda_jv_j$ with $\langle v_j\rangle$ the points of $H$, Lemma~\ref{lem:nucleus} with $Q=Q_H$ gives $\lambda_N=0$ for the column representing $N_H$, which contradicts Lemma~\ref{lem:mds}.
\end{proof}

For $q=4$, with $\F_4=\{0,1,\omega,\omega^2\}$, the only hyperoval through the frame is
\[
\{(1,0,0),(1,1,1),(1,\omega,\omega^2),(1,\omega^2,\omega),(0,1,0),(0,0,1)\},
\]
and $Q_H=x_1^2+x_2x_3$ with $N_H=(1,0,0)$. For $q=8$, with $p=x^3+x+1$, there are exactly ten hyperovals through the frame; Table~\ref{tab:q8} lists them with suitable $Q_H$ and $N_H$. Both lists were obtained by exhaustive computer search (Section~\ref{sec:comp}); the forms in Table~\ref{tab:q8} were checked separately by evaluating them at the $q+2$ points. So Proposition~\ref{prop:nonexist} applies for $q=4$ and $q=8$.

For $q\ge4$, a quadratic form that vanishes at exactly $q+1$ points of a hyperoval defines a non-degenerate conic, because a line contains at most two points of an arc and a pair of lines at most four. So the hypothesis of Proposition~\ref{prop:nonexist} says that every hyperoval through the frame is regular, and the enumeration above is a computer re-verification, for $q=4$ and $q=8$, of Segre's theorem that every hyperoval of $\PG(2,q)$ with $q\le8$ is regular~\cite{Se57}. Conversely, Segre's theorem and Corollary~\ref{cor:regular} give the non-existence for $q\le8$ without a computer. Segre's theorem also explains the number ten. For $q=8$ a regular hyperoval contains only one conic, since two distinct conics share at most four points. The conics through the frame form a pencil of $q+1$ conics, three of which are pairs of lines. So a hyperoval through the frame is either one of the other $q-2=6$ conics together with its nucleus, or the unique conic through three frame points whose nucleus is the fourth, together with that nucleus ($4$ cases; for the nucleus $(1,0,0)$ the conic is $x_1^2+x_2x_3=0$). In Table~\ref{tab:q8}, rows 4, 5, 8 and 10 are of the second kind.

\begin{table}[ht]
\caption{The ten hyperovals of $\PG(2,8)$ through the standard frame: the frame and the six points $(1,a,b)$ listed. $Q_H$ vanishes on $H\setminus\{N_H\}$ and not at $N_H$. All entries are in the integer encoding of Section~\ref{sec:1664} with $p=x^3+x+1$ (so $3\,x_1x_2$ means $(\eta+1)x_1x_2$).}\label{tab:q8}
\small
\setlength{\tabcolsep}{4pt}
\begin{tabular}{rlcl}
\hline
\# & points $(1,a,b)$ besides the frame, as $(a,b)$ & $N_H$ & $Q_H$\\
\hline
\TableQEightRows
\hline
\end{tabular}
\end{table}

\begin{proof}[Proof of Theorem~\ref{thm:main}]
For $q\in\{4,8\}$ no code exists by Proposition~\ref{prop:nonexist} and the enumeration above (alternatively, by Segre's theorem~\cite{Se57} and Corollary~\ref{cor:regular}). For $q=16$ and $q=64$ codes exist by Proposition~\ref{prop:1664}. For $q=2^m$ with $m=5$ or $m\ge7$, codes exist by Theorem~\ref{thm:family}.
\end{proof}

\section{Odd $q$: Conjecture~2}\label{sec:odd}

In this section $q$ is odd and $1\le k\le(q-1)/2$. We write $a_1,\dots,a_q$ for the elements of $\F_q$, and put $e(t)=(1,t,\dots,t^{k-1})^{T}$ for $t\in\F_q$ and $e(\infty)=(0,\dots,0,1)^{T}$. As in Lemma~\ref{lem:mds}, a $k\times n$ matrix generates an $[n,k,n-k+1]$ MDS code if and only if every $k$ of its columns are linearly independent. For $k\ge2$ the columns then represent an \emph{arc} of $\PG(k-1,q)$: $n$ distinct points, no $k$ of which lie in a hyperplane. The $q+1$ points $\langle e(t)\rangle$, $t\in\F_q\cup\{\infty\}$, form an arc. This arc and its images under $GL(k,q)$ are the \emph{normal rational curves} of $\PG(k-1,q)$ \cite{BL}. For $v\in(\F_q^*)^{q+1}$, the \emph{doubly extended generalized Reed--Solomon code} $\mathrm{GRS}_k(v)$ is the row space of $[v_1e(a_1)|\cdots|v_qe(a_q)|v_{q+1}e(\infty)]$, that is,
\[
\mathrm{GRS}_k(v)=\Bigl\{\bigl(v_1f(a_1),\dots,v_qf(a_q),v_{q+1}f_{k-1}\bigr) : f=\textstyle\sum_{j=0}^{k-1}f_jx^j\in\F_q[x]\Bigr\}.
\]

\begin{lemma}\label{lem:grsso}
$\mathrm{GRS}_k(v)$ is self-orthogonal if and only if there is a polynomial $W\in\F_q[x]$ of degree $q+1-2k$ with leading coefficient $v_{q+1}^2$ such that $W(a_i)=v_i^2$ for $1\le i\le q$.
\end{lemma}

Lemma~\ref{lem:grsso} is a reformulation of \cite[Lemma~III.2]{CL} (see also \cite[Lemma~2]{FFLZ}); we include the short proof for completeness.

\begin{proof}
Let $W=\sum_{l=0}^{q-1}W_lx^l$ be the polynomial of degree at most $q-1$ with $W(a_i)=v_i^2$ for all $i$. A polynomial as in the statement has degree $q+1-2k<q$, so it can only be this $W$. For $f,g\in\F_q[x]$ of degree less than $k$, the coefficient of $x^{2k-2}$ in $fg$ is $f_{k-1}g_{k-1}$. So the inner product of the codewords given by $f$ and $g$ is $L(fg)$, where
\[
L(h)=\sum_{i=1}^{q}W(a_i)h(a_i)+v_{q+1}^2\,[x^{2k-2}]h
\]
and $[x^{2k-2}]h$ denotes the coefficient of $x^{2k-2}$ in $h$. The products $fg$ span the polynomials of degree at most $2k-2$, so $\mathrm{GRS}_k(v)$ is self-orthogonal if and only if $L(x^j)=0$ for $0\le j\le 2k-2$. For $0\le s\le 2q-4$, $\sum_{a\in\F_q}a^s$ (with $0^0=1$) is $-1$ if $s=q-1$ and $0$ otherwise. Indeed, the proof of Lemma~\ref{lem:powersum} gives $\sum_{a\in\F_q^*}a^s=q-1=-1$ if $(q-1)\mid s$ and $0$ otherwise, and for $s=0$ the term $0^0=1$ cancels this $-1$. Since $j\le2k-2\le q-3$, this gives $\sum_iW(a_i)a_i^j=\sum_lW_l\sum_ia_i^{l+j}=-W_{q-1-j}$, and so
\[
L(x^j)=-W_{q-1-j}+v_{q+1}^2\,[j=2k-2].
\]
Hence $\mathrm{GRS}_k(v)$ is self-orthogonal if and only if $W_l=0$ for $q+2-2k\le l\le q-1$ and $W_{q+1-2k}=v_{q+1}^2$. That is, $W$ has degree $q+1-2k$ and leading coefficient $v_{q+1}^2$.
\end{proof}

\begin{proposition}\label{prop:grsodd}
If $k=(q-1)/2$, then $\mathrm{GRS}_k(v)$ is not self-orthogonal for any $v\in(\F_q^*)^{q+1}$.
\end{proposition}

\begin{proof}
Suppose that it is. By Lemma~\ref{lem:grsso} there are $A,B,C\in\F_q$ with $A=v_{q+1}^2\ne0$ such that $W(x)=Ax^2+Bx+C$ is a nonzero square for every $x\in\F_q$. Put $\delta=(B^2-4AC)/(4A^2)$. Then $W(x)=A\bigl((x+B/(2A))^2-\delta\bigr)$. As $A$ is a nonzero square, $y^2-\delta$ is a nonzero square for every $y\in\F_q$. For $\delta=0$ this fails at $y=0$. For $\delta\ne0$, every $y\in\F_q$ gives exactly two $z\in\F_q$ with $z^2=y^2-\delta$, so the equation $y^2-z^2=\delta$ has $2q$ solutions $(y,z)\in\F_q^2$. But $(y,z)\mapsto(y-z,y+z)$ is a bijection of $\F_q^2$, and $uu'=\delta$ has exactly $q-1$ solutions $(u,u')\in\F_q^2$. This is a contradiction.
\end{proof}

\begin{remark}\label{rem:grsodd}
For $1\le k\le(q-3)/2$, Lemma~\ref{lem:grsso} does give self-orthogonal codes. Let $U\in\F_q[x]$ be monic and irreducible of degree $(q+1-2k)/2\ge2$, so that $U$ has no root in $\F_q$. Then $W=U^2$ shows that $\mathrm{GRS}_k(v)$ with $v_i=U(a_i)$ and $v_{q+1}=1$ is self-orthogonal. So, among the dimensions $k\le(q-1)/2$, $k=(q-1)/2$ is the only one for which no doubly extended generalized Reed--Solomon code of length $q+1$ is self-orthogonal. This is in line with \cite[Theorem~4.8]{Liu}.
\end{remark}

\begin{proof}[Proof of Theorem~\ref{thm:odd}]
The first statement is Proposition~\ref{prop:grsodd}. Let $k=(q-1)/2$, and suppose that $C$ is a $q$-ary $[q+1,k]$ self-orthogonal MDS code with generator matrix $G=[g_1|\cdots|g_{q+1}]$. Suppose there are $A\in GL(k,q)$, a bijection $s\mapsto t_s$ from $\{1,\dots,q+1\}$ onto $\F_q\cup\{\infty\}$ and $\lambda_s\in\F_q^*$ with $Ag_s=\lambda_se(t_s)$ for all $s$. Then $C$ is the row space of $AG$, so after reordering the coordinates, $C$ becomes a doubly extended generalized Reed--Solomon code. Reordering the coordinates keeps $C$ self-orthogonal, which contradicts Proposition~\ref{prop:grsodd}.

It remains to find such $A$, $t_s$ and $\lambda_s$. For $k\ge2$ they exist exactly when the points $\langle g_s\rangle$ form a normal rational curve; this gives the last statement of the theorem. If $k=1$ ($q=3$), take $A=1$ and $\lambda_s=g_s$. If $k=2$ ($q=5$), the columns of $G$ represent $q+1$ distinct points of $\PG(1,q)$, that is, all of them, so we can take $A=I$. If $q=p$ is a prime and $k\ge3$, then $k\le p$ and $k\ne(q+1)/2$. By Ball's theorem~\cite{Ba} (see \cite[Theorem~36]{BL}), every arc of size $q+1$ in $\PG(k-1,q)$ is then a normal rational curve. If $q=9$, then $k=4$, and by Segre's theorem~\cite{Se55} (see \cite[Theorem~62]{BL}), every arc of size $10$ in $\PG(3,9)$ is a twisted cubic, that is, a normal rational curve.
\end{proof}

\section{Computations}\label{sec:comp}

All computations are elementary and exact. The program \texttt{som\_check.py} (Python~3, standard library only), which accompanies this note and contains $G_{16}$ and $G_{64}$, checks: (1)~Proposition~\ref{prop:1664}, that is, $GG^{T}=0$ and all $816$ and $45760$ column-minors are nonzero; (2)~Theorem~\ref{thm:family} exactly as stated, for $m\in\{5,7,8,9,10\}$, every admissible $i$, and three values of $c$ for each $(m,i)$; (3)~for $q=4$ and $q=8$, the exhaustive enumeration of the hyperovals through the standard frame ($1$ and $10$ of them), a form $Q_H$ as in Proposition~\ref{prop:nonexist} for each of them, and, independently, that the solution space of the six equations of Lemma~\ref{lem:search} contains no vector without zero coordinates; (4)~for every odd $q\le31$ and $k=(q-1)/2$, that the solutions $w\in\F_q^{q+1}$ of $\sum_{s}w_s\,e(t_s)e(t_s)^{T}=0$, where $t_s$ runs over $\F_q\cup\{\infty\}$, form the $3$-dimensional space $\{(W(a_1),\dots,W(a_q),W_2):\deg W\le2\}$ given by Lemma~\ref{lem:grsso}, and that none of them has all coordinates nonzero squares (Proposition~\ref{prop:grsodd}); as a control, the codes of Remark~\ref{rem:grsodd} with $k=(q-3)/2$ are self-orthogonal. The same results were also confirmed by separately written programs.

\subsection*{Acknowledgment}
This work was developed with the assistance of Anthropic's Claude Opus~5.5. The author takes full responsibility for the content, verification and correctness of all results presented here.

\end{document}

%% file: tables.tex
\newcommand{\MatrixSixteen}{%
\begin{pmatrix}
14 & 12 & 3 & 12 & 10 & 1 & 11 & 2 & 2 & 7 & 1 & 11 & 5 & 5 & 11 & 11 & 0 & 0 \\
0 & 12 & 6 & 7 & 14 & 5 & 15 & 14 & 3 & 10 & 10 & 9 & 9 & 12 & 8 & 3 & 3 & 0 \\
0 & 12 & 15 & 5 & 9 & 7 & 5 & 6 & 5 & 3 & 12 & 6 & 11 & 14 & 3 & 8 & 0 & 9
\end{pmatrix}}

\newcommand{\MatrixSixtyFourA}{%
\begin{pmatrix}
46 & 43 & 54 & 8 & 13 & 60 & 31 & 42 & 22 & 56 & 17 & 33 & 39 & 56 & 9 & 55 & 40 & 38 & 9 & 6 & 23 & 55 \\
0 & 43 & 47 & 24 & 52 & 9 & 1 & 19 & 54 & 49 & 44 & 36 & 46 & 20 & 61 & 49 & 30 & 29 & 4 & 41 & 32 & 9 \\
0 & 43 & 42 & 47 & 25 & 27 & 30 & 49 & 41 & 54 & 21 & 52 & 20 & 52 & 2 & 24 & 35 & 17 & 35 & 24 & 56 & 10
\end{pmatrix}}
\newcommand{\MatrixSixtyFourB}{%
\begin{pmatrix}
5 & 62 & 43 & 30 & 45 & 16 & 32 & 52 & 57 & 29 & 18 & 3 & 27 & 5 & 7 & 7 & 30 & 45 & 63 & 36 & 60 & 62 \\
13 & 13 & 57 & 2 & 51 & 58 & 18 & 28 & 19 & 36 & 27 & 32 & 1 & 41 & 57 & 62 & 22 & 28 & 48 & 62 & 48 & 25 \\
38 & 30 & 53 & 51 & 32 & 37 & 59 & 14 & 20 & 57 & 41 & 30 & 14 & 31 & 49 & 16 & 5 & 51 & 37 & 20 & 20 & 1
\end{pmatrix}}
\newcommand{\MatrixSixtyFourC}{%
\begin{pmatrix}
62 & 25 & 39 & 29 & 57 & 5 & 17 & 55 & 3 & 34 & 19 & 32 & 13 & 4 & 62 & 25 & 32 & 19 & 20 & 27 & 0 & 0 \\
37 & 7 & 54 & 15 & 6 & 48 & 6 & 48 & 31 & 39 & 46 & 13 & 20 & 33 & 22 & 8 & 34 & 5 & 12 & 18 & 45 & 0 \\
54 & 55 & 3 & 49 & 2 & 6 & 1 & 21 & 59 & 17 & 33 & 12 & 26 & 26 & 13 & 62 & 28 & 12 & 25 & 41 & 0 & 61
\end{pmatrix}}

\newcommand{\TableQEightRows}{%
1 & (2,3), (3,5), (4,2), (5,6), (6,7), (7,4) & $(1,2,3)$ & $3\,x_1x_2 + 2\,x_1x_3 + x_2x_3$ \\
2 & (2,3), (3,6), (4,5), (5,7), (6,4), (7,2) & $(1,4,5)$ & $5\,x_1x_2 + 4\,x_1x_3 + x_2x_3$ \\
3 & (2,4), (3,2), (4,7), (5,3), (6,5), (7,6) & $(1,3,2)$ & $2\,x_1x_2 + 3\,x_1x_3 + x_2x_3$ \\
4 & (2,4), (3,5), (4,6), (5,7), (6,2), (7,3) & $(0,1,0)$ & $x_2^2 + x_1x_3$ \\
5 & (2,5), (3,6), (4,7), (5,2), (6,3), (7,4) & $(1,0,0)$ & $x_1^2 + x_2x_3$ \\
6 & (2,5), (3,7), (4,3), (5,4), (6,2), (7,6) & $(1,7,6)$ & $6\,x_1x_2 + 7\,x_1x_3 + x_2x_3$ \\
7 & (2,6), (3,4), (4,5), (5,2), (6,7), (7,3) & $(1,6,7)$ & $7\,x_1x_2 + 6\,x_1x_3 + x_2x_3$ \\
8 & (2,6), (3,7), (4,2), (5,3), (6,4), (7,5) & $(0,0,1)$ & $x_3^2 + x_1x_2$ \\
9 & (2,7), (3,2), (4,6), (5,4), (6,3), (7,5) & $(1,5,4)$ & $4\,x_1x_2 + 5\,x_1x_3 + x_2x_3$ \\
10 & (2,7), (3,4), (4,3), (5,6), (6,5), (7,2) & $(1,1,1)$ & $x_1x_2 + x_1x_3 + x_2x_3$ \\}

%% file: main.bbl
\begin{thebibliography}{CPPR}

\bibitem[Ba]{Ba}
S.~Ball, \emph{On sets of vectors of a finite vector space in which every subset of basis size is a basis}, J.\ Eur.\ Math.\ Soc.\ \textbf{14} (2012), 733--748.

\bibitem[BL]{BL}
S.~Ball and M.~Lavrauw, \emph{Arcs in finite projective spaces}, EMS Surv.\ Math.\ Sci.\ \textbf{6} (2019), 133--172; theorem numbers refer to \href{https://arxiv.org/abs/1908.10772v2}{arXiv:1908.10772v2}.

\bibitem[CL]{CL}
B.~Chen and H.~Liu, \emph{New constructions of MDS codes with complementary duals}, IEEE Trans.\ Inform.\ Theory \textbf{64} (2018), 5776--5782.

\bibitem[CPPR]{CPPR}
W.~Cherowitzo, T.~Penttila, I.~Pinneri and G.~F.~Royle, \emph{Flocks and ovals}, Geom.\ Dedicata \textbf{60} (1996), 17--37.

\bibitem[FFLZ]{FFLZ}
W.~Fang, F.-W.~Fu, L.~Li and S.~Zhu, \emph{Euclidean and Hermitian hulls of MDS codes and their applications to EAQECCs}, IEEE Trans.\ Inform.\ Theory \textbf{66} (2020), 3527--3537.

\bibitem[Hir]{Hir}
J.~W.~P.~Hirschfeld, \emph{Projective Geometries over Finite Fields}, 2nd ed., Oxford Mathematical Monographs, Clarendon Press, Oxford, 1998.

\bibitem[Liu]{Liu}
S.~Liu, \emph{On Self-orthogonal MDS Codes}, \href{https://arxiv.org/abs/2611.30048v1}{arXiv:2611.30048v1} [cs.IT], 2026.

\bibitem[LS]{LS}
L.~Lunelli and M.~Sce, \emph{$k$-archi completi nei piani proiettivi desarguesiani di rango $8$ e $16$}, Centro di Calcoli Numerici, Politecnico di Milano, 1958.

\bibitem[MJ]{MJ}
E.-J.~Moreno-Mej\'{\i}a and H.~Jim\'enez-Hern\'andez, \emph{Code and data for ``Self-orthogonal MDS codes of length $q+2$ and $q+1$''}, Zenodo, 25 November 2026, \href{https://doi.org/10.5281/zenodo.24526969}{doi:10.5281/zenodo.24526969}.

\bibitem[Se55]{Se55}
B.~Segre, \emph{Curve razionali normali e $k$-archi negli spazi finiti}, Ann.\ Mat.\ Pura Appl.\ (4) \textbf{39} (1955), 357--379.

\bibitem[Se57]{Se57}
B.~Segre, \emph{Sui $k$-archi nei piani finiti di caratteristica due}, Rev.\ Math.\ Pures Appl.\ \textbf{2} (1957), 289--300.

\end{thebibliography}
